\documentclass[preprint, 12pt, 3p]{elsarticle}

\usepackage{amsmath, amssymb, amsfonts, amsthm}
\RequirePackage[numbers]{natbib}
\RequirePackage[colorlinks,citecolor=blue,urlcolor=blue]{hyperref}

\theoremstyle{plain}
\newtheorem{theorem}{Theorem}[section]
\newtheorem{lemma}[theorem]{Lemma}
\newtheorem{assumption}{Assumption}[section]
\newtheorem{corollary}[theorem]{Corollary}
\newtheorem{remark}{Remark}[section]
\newtheorem{definition}{Definition}[section]

\def\B{\mathfrak{B}}
\def\C{\mathcal}
\def\BB{\mathbb}

\makeatletter
\def\ps@pprintTitle{%
  \let\@oddhead\@empty
  \let\@evenhead\@empty
  \let\@oddfoot\@empty
  \let\@evenfoot\@oddfoot
}
\makeatother

\allowdisplaybreaks

\journal{Journal of Mathematical Analysis and Applications}

\begin{document}

\begin{frontmatter}

\title{On solutions of Kolmogorov's equations for nonhomogeneous jump Markov processes}

\author[SUNY]{Eugene A. Feinberg\fnref{FM}}
\ead{eugene.feinberg@stonybrook.edu}

\author[SUNY]{Manasa Mandava\fnref{FM}}
\ead{mmandava@ams.sunsyb.edu}

\fntext[FM]{This research was partially supported by NSF grant CMMI-0928490}
\address[SUNY]{Department of Applied Mathematics and Statistics, Stony Brook University, Stony Brook, USA 11794-3600}

\author[STEKLOV]{Albert N. Shiryaev}
\ead{albertsh@mi.ras.ru}

\address[STEKLOV]{Steklov Mathematical Institute, 8, Gubkina Str., Moscow, Russia 119991}

\begin{abstract}
This paper studies three ways to construct a nonhomogeneous jump
Markov process: (i) via a compensator of the random measure of a
multivariate point process, (ii) as a minimal solution of the
backward Kolmogorov equation, and (iii) as a minimal solution of
the forward Kolmogorov equation.  The main conclusion of this
paper is that, for a given measurable transition intensity,
commonly called a $Q$-function, all these constructions define the
same transition function.  If this transition function is regular,
that is, the probability of accumulation of jumps is zero, then
this transition function is the unique solution of the backward
and forward Kolmogorov equations. For continuous $Q$-functions,
Kolmogorov equations were studied in Feller's seminal paper. In
particular, this paper extends Feller's results for continuous
$Q$-functions to measurable $Q$-functions and provides additional
results.
\end{abstract}

%

\end{frontmatter}

\section{Introduction}\label{sec1}

Let $(X, \mathfrak{B}(X))$ be a standard Borel space, that is,
$(X, \mathfrak{B}(X))$ is a measurable space for which there
exists a measurable injection onto a Borel subset of the  real
line  endowed with its Borel $\sigma$-field.  For a Borel subset
$E$ of the extended real line, we denote by $\mathfrak{B}(E)$ its
Borel $\sigma$-field.
A function $P(u, x; t, B)$, where $u,t\in \BB{R}_+:= ]0,
\infty[   ,$ $u<t$, $x\in X$, and $B\in\B(X)$, is called  a
transition function if it takes values in $[0,1]$
and satisfies the following properties:
\begin{itemize}
\item[(i)] for all $u,x,t$ the function $P(u,x;t,\cdot)$ is a
measure on $(X, \mathfrak{B}(X))$;
 \item[(ii)] for all $B$ the
function $P(u,x; t, B)$ is Borel measurable in $(u,x,t);$ \item[(iii)]
$P(u,x;t,B)$ satisfies the Chapman-Kolmogorov equation
\begin{equation}
\label{CKE} P(u, x; t, B) = \int_{X} P(s,y; t, B)P(u, x; s, dy),
\qquad u < s < t.
\end{equation}
\end{itemize}
If $P(u,x;t,X) = 1$ for all $u,x,t$, then the transition function $P$ is called a {\it
regular transition function}.

A stochastic process $\{\BB{X}_t: t \ge 0\}$ with values in $X$,
defined on a probability space $(\Omega, \C{F}, \BB{P})$ with a
filtration $\{\C{F}_t\}_{t \ge 0}$, is called Markov if
$\BB{P}(\BB{X}_t \in B \mid \C{F}_u) = \BB{P}(\BB{X}_t \in B \mid
\BB{X}_u)$, $\BB{P}-a.s$ for all $u,t \in \BB{R}_+$ with $u <t$ and for all $B
\in \B(X)$. And each Markov process has a transition function $P$ such
that $\BB{P}(\BB{X}_t \in B \mid \BB{X}_u) = P(u,\BB{X}_u; t, B),$
$\BB{P}-a.s$; 
see Kuznetsov \cite{Kuz}, where the equivalence of the two
approaches to define a Markov process from Kolmogorov \cite{Kol}
is established. If each sample path of the Markov process is a
right continuous piecewise constant function that has a finite
number of discontinuity points on each interval $[0,t]$ for
$t<t_\infty$, where $t_\infty$ is the lifetime of the process,
then the Markov process is called a {\it jump Markov process};
Gikhman and Skorokhod \cite[p. 187]{Gik}. In this paper,
measurability and Borel measurability are used synonymously and
all conditional probabilities are defined almost sure, even when
this is not explicitly stated.

We now introduce the following function $q$ that can be
interpreted as the transition intensities of a nonhomogeneous jump
Markov process.
A function $q(x,t,B)$, where $x \in X$, $t \in \BB{R}_+$, and $B \in \B(X)$,
is called a {\it Q-function} if it satisfies the following
properties:
\begin{itemize}
\item[(i)]for all $x,t$ the function $q(x,t,\cdot)$ is a signed
measure on $(X, \mathfrak{B}(X))$ such that $q(x,t,X)$ $\leq$ $0$ and
$0 \leq q(x, t, B \setminus \{x\}) < \infty$ for all $B \in
\mathfrak{B}(X)$;
\item[(ii)] for all $B$ the
function $q(x,t,B)$ is measurable in $(x,t).$
\end{itemize}
In addition to properties (i) and (ii), if $q(x,t,X) =0$ for all
$x,t$, then the $Q$-function $q$ is  called {\it conservative}.
Note that any $Q$-function can be transformed into a conservative
$Q$-function  by adding a state $\tilde{x}$ to $X$ with
$q(x,t,\{\tilde{x}\}):= -q(x,t,X)$, $q(\tilde{x}, t, X) := 0$, and
$q(\tilde{x}, t, \{\tilde{x}\}) := 0$, where $x \in X$ and $t \in
\mathbb{R}_+$. To simplify the presentation, in this paper we
always assume that $q$ is conservative. If there is no assumption
that $q$ is conservative, Remark~\ref{rm4.1} explains how the main
formulations change. A $Q$-function $q$ is called {\it continuous}
if it is continuous in $t \in \BB{R}_+$.

Let $q(x,t) := - q(x,t,\{x\})$. A set $B \in \mathfrak{B}(X)$ is called {\it q-bounded} if $\sup_{x \in B, t \in \mathbb{R}_+}q(x,t) < \infty$ and the $Q$-function $q$ is called {\it stable} if the set $\{x\}$ is $q$-bounded for all $x \in X$. The following assumption introduced by Feller~\cite{Fel} holds throughout the paper.
\begin{assumption}
\label{A}
The $Q$-function $q$ is stable.
\end{assumption}

Let $B_n := \{x \in X: \sup_{t \in \mathbb{R}_+}q(x,t) < n+1\}$
for all $n\geq 0$. If the $Q$-function $q$ is stable, then $B_n
\uparrow X$ as $n \to \infty$. Thus, Assumption \ref{A} is
equivalent to the existence of a sequence of $q$-bounded sets
$\{B_n\}_{n \geq 0}$ such that $B_n \uparrow X$ as $n \to \infty$.
This way a stable $Q$-function was defined in Feller \cite{Fel}.
In this paper, a non-negative solution $\bar{f}$ in a certain
class of solutions of a functional equation is called the minimal
non-negative solution if for any non-negative solution $f$ of this
equation from that class $\bar{f}(\textbf{x}) \le f(\textbf{x})$
for all values of the argument $\textbf{x}$.

Feller~\cite{Fel} studied the backward and forward Kolmogorov
equations for continuous $Q$-functions.    For a stable continuous
$Q$-function, Feller~\cite{Fel} provided explicit formulae for a
transition function that satisfies both the backward and forward
Kolmogorov equations. If the constructed transition function is
regular, Feller~\cite[Theorem 3]{Fel} showed that this transition
function is the unique non-negative solution to the backward
Kolmogorov equation. Though Feller~\cite{Fel} focused on regular
transition functions,  it follows from the proof of Theorem 3 in
Feller~\cite{Fel} that the  transition function constructed there
is the minimal
non-negative solution to the backward Kolmogorov equation. 
For homogeneous Markov processes, that is $Q$-functions do not
depend on the time parameter, Doob~\cite{DoobP}, \cite[Chap.
6]{Doob} provided an explicit construction for multiple transition
functions satisfying the backward Kolmogorov equation. For
countable-state homogeneous Markov processes, Kendall~\cite{Ken},
Kendall and Reuter~\cite{KenReu}, and Reuter~\cite{Reu} gave
examples with non-unique solutions to Kolmogorov equations and
Reuter~\cite{Reu} provided necessary and sufficient conditions for
their uniqueness; see also Anderson~\cite{And} and Chen et
al.~\cite{Che}. Ye, Guo, and Hern\'{a}ndez-Lerma~\cite{Guo} proved
the existence of a transition function that is the minimal
non-negative solution to both the backward and forward Kolmogorov
equations for a countable state problem with measurable
$Q$-functions. A conservative $Q$-function can be used to
construct a predictable random measure. According to
Jacod~\cite[Theorem 3.6]{Jac}, an initial state distribution and a
predictable random measure define uniquely a  multivariate point
process.

This paper studies the backward and forward Kolmogorov equations
for measurable $Q$-functions and standard Borel state spaces. It
extends Feller's~\cite{Fel} results for continuous $Q$ functions
to measurable $Q$-functions and establishes additional results.
Theorem~\ref{JMP} below states that the stochastic process
associated with the multivariate point process defined by a stable
conservative $Q$ function  and an initial state distribution is a
jump Markov process with the transition function defined in
Feller~\cite{Fel}.
 Theorems~\ref{BKE} and \ref{FKE} state
that this transition function  satisfies the backward and forward
Kolmogorov equations. In addition, this transition function is the
minimal non-negative solution of both the backward and forward
Kolmogorov equations and, if this transition function is regular,
then it is the unique non-negative solution of the backward and
forward Kolmogorov equations; Theorems~\ref{unique-b},
\ref{unique-f}. Thus, the minimal non-negative solution of both
the backward and forward Kolmogorov equations  is the  transition
function of a jump Markov process associated with a multivariate
point process whose compensator is defined by the conservative
$Q$-function.

In addition to answering the fundamental question on how to
construct a jump Markov process with a given $Q$-function, our
interest in this study is motivated by its applications to control
of continuous-time jump Markov  processes. Here we mention two of
them:

 (i) For a countable state space,
 each Markov policy along with a given initial state
distribution defines a jump Markov process with the transition
function being the minimal non-negative solution of the forward
Kolmogorov equation; Guo and Hern\'andez-Lerma~\cite[Section
2.2]{Her}. An arbitrary policy defines a multivariate point
process via the compensator of its random measure;
Kitaev~\cite{Kit},
 Kitaev and Rykov~\cite[Section
4.6]{KR}, Feinberg~\cite{Fei, Fei1}, Guo and Piunovskiy~\cite{GP}.
The results of this paper imply that for Markov policies these two
definitions are equivalent for problems with Borel state spaces.

(ii) Feller's~\cite{Fel} results are broadly used in the
literature on  continuous-time Markov decision processes to define
Markov processes corresponding to Markov policies, and this leads
to the unnecessary assumption that decisions depend continuously
on time; see, e.g., Guo and Rieder \cite[Definition 2.2]{GuoR}.
For countable state problems, the results of Ye, Guo, and
Hern\'{a}ndez-Lerma~\cite{Guo} removed the necessity to assume
this continuity.  The results of the current paper imply that this
continuity assumption is unnecessary for Markov decision processes
with Borel state spaces.

\section{Relation between Jump Markov Processes and Q-functions}
\label{Q-J}

The main goals of this section are to show that an
initial state distribution and a (stable) $Q$-function $q$
define a jump Markov process and to construct
its transition function.

Let $x_\infty\notin X$ be an isolated point adjoined to the space
$X$. Denote ${\bar X}=X\cup\{x_\infty\}$ and $\bar
{\BB{R}}_+=\ ]0,\infty]$. Consider the Borel $\sigma$-field $\B({\bar
X})=\sigma(\mathfrak{B}(X),\{x_\infty\})$ on $\bar X$, which is
the minimal $\sigma$-field containing $ \mathfrak{B}(X)$ and
$\{x_\infty\}.$ Let $({\bar X} \times \bar
{\BB{R}}_+)^\infty$ be the set of all sequences $(x_0, t_1, x_1,
t_2, x_2, \ldots)$ with $x_n\in \bar{X}$ and $t_{n+1}\in
\bar{\BB{R}}_+$ for all $n \geq 0.$ This set is endowed
with the $\sigma$-field generated by the products of the Borel
$\sigma$-fields $\B(\bar{X})$ and $\B(\bar{\BB{R}}_+)$.

Denote by $\Omega$ the subset of all sequences $\omega= (x_0, t_1,
x_1, t_2, x_2, \ldots)$ from $({\bar X} \times \bar
{\BB{R}}_+)^\infty$ such that: (i) $x_0 \in X$; (ii) if $t_n < \infty$, then $t_n < t_{n+1}$ and
$x_n \in X$, and, if $t_n = \infty$, then $t_{n+1} = t_n$ and
$x_n = x_\infty$, for all $n\geq 1$. Observe that $\Omega$ is a measurable subset of
$({\bar X} \times \bar
{\BB{R}}_+)^\infty$. Consider the measurable space $(\Omega,
\C{F})$, where $\C{F}$ is the $\sigma$-field of the measurable subsets
of $\Omega$. Then, $x_n(\omega)=x_n$ and
$t_{n+1}(\omega)=t_{n+1}$, $n \ge 0$, are random variables defined on the
measurable space $(\Omega, \C{F})$. Let $t_0 := 0$, $t_\infty
(\omega) := \lim\limits_{n \to \infty} t_n (\omega)$, $\omega \in \Omega$, and for all $t \ge 0$ let $\C{F}_t := \sigma(\B(X), \C{G}_t)$, where $\C{G}_t := \sigma (I\{x_n \in B\}I\{t_n \le s\}: n \ge 1, 0 \le s \le t, B \in \B(X)).$ Throughout this paper, we omit $\omega$ whenever possible and also follow the standard convention that $0
\times \infty = 0$.

For a given $Q$-function $q$, consider the random measure
$\nu$ on $(\BB{R}_+ \times X, \B(\BB{R}_+) \times \B(X))$ defined
by
\begin{equation}
\label{compensator}
\nu(\omega; ]0,t], B) =
\int_{0}^{t}\sum\limits_{n \ge 0}I\{t_n < s \leq t_{n+1}\}q(x_n,
s, B \setminus \{x_n\}) ds,\quad  t \in \BB{R}_+, \ B \in
\B(X).
\end{equation}
Note that $\nu(\{t\} , X)= \nu([t_\infty, \infty[\, ,
X)= 0$ and \eqref{compensator} can be rearranged as
\begin{multline}
\label{e21}
\nu(]0,t] , B)= \sum_{n \ge 0} I\{t_n < t \le t_{n+1} \} \left (\sum\limits_{m=0}^{n-1} \int_0^{t_{m+1} - t_m} q(x_m, t_m+s, B \setminus \{x_m\})ds \right . \\
\left. + \int_0^{t-t_n} q(x_n, t_n+s, B \setminus \{x_n\})ds\right).
\end{multline}
As the expression in the parentheses on the right hand side of
\eqref{e21} is an $\C{F}_{t_n}$-measurable process for each
$B\in\B(X)$, it follows from Jacod \cite[Lemma 3.3]{Jac}  that the
process $\{\nu(]0,t], B): t \in \BB{R}_+\}$ is predictable.
Therefore, the measure $\nu$ is a predictable random measure.
According to Jacod \cite[Theorem 3.6]{Jac}, the predictable random
measure $\nu$ defined in \eqref{compensator} and a probability
measure $\mu$ on $X$ define a unique probability measure $\BB{P}$
on $(\Omega, \C{F})$ such that $\BB{P}(x_0 \in B) = \mu(B)$, $B \in \B(X)$,
and $\nu$ is the compensator of the random
measure of the multivariate point process $(t_n, x_n)_{n \ge 1}$
defined by the triplet $(\Omega,\C{F},\BB{P})$.

Consider the process $\{\BB{X}_t: t \ge 0\}$,
\begin{equation}
\label{jump}
\mathbb{X}_t(\omega) := \sum_{n \geq 0} I\{t_n \leq t < t_{n+1}\}x_n  + I\{t_\infty \leq
t\}x_\infty,
\end{equation}
defined on $(\Omega, \C{F}, \BB{P})$. We abbreviate the process $\{\BB{X}_t: t\ge 0\}$ as $\BB{X}$. The main result of this section, Theorem \ref{JMP}, shows that the process $\BB{X}$ is a jump Markov process and provides its transition function.

For an $\C{F}_t$-measurable stopping time $\tau$, let $N(\tau) := \max\{n = 0,1,\ldots: \tau \ge t_n\}$. Since  $N(\tau) = \infty$ and $\BB{X}_\tau = \{x_\infty\}$ when $\tau \ge t_\infty$, we follow the convention that $t_{\infty+1} = \infty$ and $x_{\infty+1} = x_\infty$. Denote by
$G_{\tau}(\omega; \cdot, \cdot)$ and $H_{\tau}(\omega; \cdot)$ respectively the
regular conditional laws of $(t_{N(\tau)+1}, x_{N(\tau)+1})$ and $t_{N(\tau)+1}$ with respect to $\C{F}_{\tau}$; $H_{\tau}(\omega; \cdot) = G_{\tau}(\omega; \cdot, \bar{X})$. In particular,  $G_{t_n}(\omega; \cdot, \cdot)$ and $H_{t_n}(\omega; \cdot)$, where $n = 0,1,\ldots,$ denote the conditional laws of $(t_{n+1}, x_{n+1})$ and $t_{n+1}$ with respect to $\C{F}_{t_n}$. We remark that the notations $G_{t_n}$ and $H_{t_n}$ correspond to the notations $G_n$ and $H_n$ in Jacod \cite[p. 241]{Jac}.

\begin{lemma}
For all $u,t \in \BB{R}_+$, $u < t$,
\begin{align}
\label{con4a}
H_{u}([t, \infty]) &= e^{-\int_u^t q(\BB{X}_u, s) ds}, \quad &N&(u) < \infty, \\
\label{con3a}
G_{u}(dt, B) &= e^{-\int_{u}^t q(\BB{X}_u, s)ds} q(\BB{X}_u, t, B \setminus \{\BB{X}_u\}) dt, \qquad  &B& \in \B(X),\, N(u) < \infty.
\end{align}
\begin{proof}
According to Jacod \cite[Proposition 3.1]{Jac},
for all $t \in \BB{R}_+$, $B \in \B(X)$, and $n = 0,1, \ldots$
\begin{equation}
\label{e22}
\nu(dt, B) = \frac{G_{t_n}(dt, B)}{H_{t_n}([t, \infty])}, \qquad t_n < t \le t_{n+1}.
\end{equation}
In particular, for $B=X$, from \eqref{e22} and from the property that $x_{n+1} \in X $ when $t_{n+1} < \infty$,
\begin{equation*}
\nu(dt, X) = \frac{G_{t_n}(dt, X)}{H_{t_n}([t, \infty])} =
\frac{H_{t_n}(dt)}{H_{t_n}([t, \infty])}, \qquad t_n < t \le t_{n+1}.
\end{equation*}
This equality implies that $\nu(dt, X)$ is the hazard rate
function corresponding to the distribution $H_{t_n}$ when $t_n < t \le t_{n+1}$. Therefore,
\begin{equation}
\label{e24}
H_{t_n}([t, \infty]) = e^{-\nu(]t_n, t], X)I\{t_n < t \le t_{n+1}\}}, \qquad t \in \BB{R}_+, t > t_n.
\end{equation}
From \eqref{compensator} and \eqref{e24}, for all $t \in \BB{R}_+$,
\begin{equation}
\label{con1a} H_{t_n}([t, \infty]) = e^{-\int_{t_n}^t q(x_n,s)ds}, \qquad \qquad t > t_n,
\end{equation}
and from \eqref{compensator}, \eqref{e22}, and \eqref{con1a}, for all $ t \in \BB{R}_+$, $B \in \B(X)$,
\begin{equation}
\label{con2a} G_{t_n}(dt, B) = e^{-\int_{t_n}^t q(x_n, s)d s} q(x_n, t, B \setminus \{x_n\})dt, \qquad t > t_n.
\end{equation}

To compute $G_u$, observe that for all $u,t \in \BB{R}_+$, $u <t$, and $B \in \B(X)$,
\begin{equation}
\label{e25}
\begin{aligned}
G_u(dt, B) &= \BB{P}(t_{N(u)+1} \in [t, t+dt [\, , x_{N(u)+1} \in B \mid \C{F}_u) \\
&= \sum\limits_{n \ge 0} \BB{P}(t_{N(u)+1} \in [t, t+dt [\, , x_{N(u)+1} \in B \mid \C{F}_{u})I\{N(u) = n\} \\
&= \sum\limits_{n \ge 0} \BB{P}(t_{n+1} \in [t, t+dt [\, , x_{n+1} \in B \mid \C{F}_{u}, N(u)=n )I\{N(u) = n\}, \\
\end{aligned}
\end{equation}
where the first equality follows from the definition of
$G_u$, the second equality holds because $\{N(u) = \infty\} \cup \{N(u)= n\}_{n =
0,1,\ldots}$ is an $\C{F}_u$-measurable partition of $\Omega$ and $x_{N(u)+1}=x_{\infty}\notin X$ when
$N(u)=\infty$, and the third equality follows from $N(u)=n$ and from $\{N(u) = n \} \in \C{F}_u$.

Observe that for any random variable $Z$ on $(\Omega, \C{F})$
\begin{multline}
\label{e26}
\BB{P}(Z \mid \C{F}_u, N(u) = n)I\{N(u)= n\}= \BB{P}(Z \mid \C{F}_{t_n}, N(u) = n)I\{N(u)= n\}\\
\begin{aligned}
&= \BB{P}(Z \mid \C{F}_{t_n}, t_n \le u, t_{n+1}> u)I\{N(u)= n\} = \BB{P}(Z \mid \C{F}_{t_n}, t_{n+1}> u)I\{N(u)= n\}\\
&= \frac{\BB{P}(Z,  t_{n+1} > u \mid \C{F}_{t_n})}{\BB{P}(t_{n+1} > u \mid \C{F}_{t_n})}I\{N(u) = n\},
\end{aligned}
\end{multline}
 where the first equality follows from Br\'emaud \cite[Theorem T32, p. 308]{Bremaud}, the second equality holds because
$\{t_n \le u, t_{n+1} > u\}= \{N(u) = n\}$, the third equality holds because $\{t_n \le u\} \in \C{F}_{t_n}$,  and the last
one follows from the definition of conditional probabilities. Let $Z = \{t_{n+1} \in [t, t+dt [\, , x_{n+1} \in B \}$, where $t \in \BB{R}_+, B \in \B(X)$. Then \eqref{e25} and \eqref{e26} imply
\begin{equation}
\label{e27}
\begin{aligned}
G_u(dt, B) &= \sum\limits_{n \ge 0} \frac{\BB{P}(t_{n+1} \in [t, t+dt [\, , x_{n+1} \in B \mid \C{F}_{t_n})}{\BB{P}(t_{n+1} > u \mid \C{F}_{t_n})}I\{N(u) = n\} \\
&= \sum\limits_{n \ge 0} \frac{ e^{-\int_{t_n}^t q(x_n, s)d s }q(x_n,t, B \setminus \{x_n\})dt}{e^{-\int_{t_n}^u q(x_n,s)ds}}I\{N(u) = n\}  \\
&= e^{-\int_{u}^t  q(\BB{X}_u, s)ds }q(\BB{X}_u,t, B \setminus \{\BB{X}_u\})dt,
\end{aligned}
\end{equation}
where the first equality holds because $\{t_{n+1} \in [t, t+dt [\, ,t_{n+1} > u\} = \{t_{n+1} \in [t, t+dt [\}$ when $t > u$,
the second equality follows from \eqref{con1a} and \eqref{con2a}, and the last equality holds since $x_n = \BB{X}_u$ when $N(u) = n$. For all $u,t \in \BB{R}_+$, $u <t$, it follows from the
property that $x_{N(u)+1} \in X$ when $t_{N(u)+ 1} < \infty$ and
from \eqref{e27} that $H_u([t,\infty])$ satisfies
\eqref{con4a}.
\end{proof}
\end{lemma}

Following Feller \cite[p. 501]{Fel}, for $x \in X$, $u,t \in \mathbb{R}_+$, $u < t$, and $B \in \B(X)$,
define
\begin{equation}
\label{b0} \bar{P}^{(0)} (u,x;t,B) = I\{x \in B\} e^{-\int_u^t
q(x, s) ds},
\end{equation}
and for $n \geq 1$ define
\begin{equation}
\label{bn}
\bar{P}^{(n)}(u, x; t, B) = \int_{u}^{t} \int_{X \setminus \{x\}}
 e^{ -\int_{u}^{s} q(x,\theta) d\theta }  q(x,s,dy) \bar{P}^{(n-1)}(s, y; t, B)
 ds.
\end{equation}
Set
\begin{equation}
\label{def}
\bar{P}(u, x; t, B) := \sum\limits_{n=0}^{\infty} \bar{P}^{(n)}(u, x; t, B).
\end{equation}
Observe that $\bar{P}$ is a transition function. For stable continuous $Q$-functions, Feller \cite[Theorems 2, 5]{Fel} proved that
(a) for fixed $u,x,t$ the function $\bar{P}(u,x;t,\cdot)$ is a measure on $(X,\B(X))$  such that $0 \le \bar{P}(u,x;t,\cdot) \le 1$, and (b) for all $u,x,t,B$ the function $\bar{P}(u,x;t,B)$  satisfies the Chapman-Kolmogorov equation \eqref{CKE}. The proofs remain correct for measurable $Q$-functions $q$. The measurability of $\bar{P}(u,x;t,B)$ in $u,x,t$ for all $B \in \B(X)$ is straightforward from the definitions \eqref{b0}, \eqref{bn}, and \eqref{def}. Therefore, the function $\bar{P}$ satisfies properties (i)-(iii) from the definition of a transition function.

\begin{theorem}
\label{JMP} For a given initial state distribution and for a
stable $Q$-function $q$, the process $\BB{X}$
defined in \eqref{jump} is a jump Markov process with the
transition function $\bar P$.
\begin{proof}
Observe that the sample paths of the process $\BB{X}$ are
right-continuous piecewise-constant functions that have finite
number of discontinuities on each interval $[0,t]$ for $t <
t_\infty$. Thus if, for all $u,t \in \BB{R}_+$, $u < t$, and $B \in \B(X)$,
\begin{equation}
\label{e:MarkProp11}
\BB{P}(\mathbb{X}_t \in B \mid \C{F}_u) = \BB{P}(\mathbb{X}_t \in B  \mid \BB{X}_u) = \bar{P}(u,\BB{X}_u; t, B), \qquad u < t_\infty,
\end{equation}
then the process $\BB{X}$ is a jump Markov process with the transition function $\bar{P}$.
To prove  \eqref{e:MarkProp11}, we first establish by induction
that for all $n=0,1,\ldots,$ $u, t \in \BB{R}_+,$ $u < t $, and $B \in
\B(X)$
\begin{equation}
\label{iii6}
\mathbb{P}(\mathbb{X}_t \in B, N_{]u,t]} = n \mid \C{F}_u) =
\bar{P}^{(n)} (u,\BB{X}_u;t,B),  \qquad u < t_\infty,
\end{equation}
where $N_{]u,t]}:= N(t) - N(u)$ when $ u < t_\infty$ and $N_{]u,t]}:= \infty$ when $u \ge t_\infty$.
Equation \eqref{iii6} holds for $n=0$ because for $u < t_\infty$
\begin{multline}
\label{i3}
\mathbb{P}(\mathbb{X}_t \in B, N_{]u,t]} = 0 \mid \C{F}_u) = \BB{P}(\BB{X}_u \in B, t_{N(u)+ 1} > t \mid \C{F}_u)  \\
= I\{\BB{X}_u \in B\} H_u(]t, \infty]) = I\{\BB{X}_u \in B\} e^{-\int_u^t q(\BB{X}_u, s) ds} = \bar{P}^{(0)} (u,\BB{X}_u;t,B),
\end{multline}
where the first equality holds because the corresponding events
coincide, the second equality holds because
$\{\BB{X}_u \in B\} \in \C{F}_u$ and from the definition of $H_u$, the third equality
is correct because of (\ref{con4a}), and the last equality is
\eqref{b0}.

For some $n \geq 0$, assume that \eqref{iii6} holds. Then for $u < t_\infty$
\begin{multline}
\label{i4}
\mathbb{P}(\mathbb{X}_t \in B, N_{]u,t]} = n+1 \mid\C{F}_u)  \\
\begin{aligned}
&= \int_u^t \int_{X \setminus \{\BB{X}_u\}} \BB{P}( \BB{X}_t \in B, N_{]t_{N(u)+1}, t]} = n  \mid \C{F}_u, t_{N(u) + 1}, x_{N(u)+1}) G_u(dt_{N(u) + 1}, dx_{N(u)+1})\\
&= \int_u^t \int_{X \setminus \{\BB{X}_u\}} \BB{P}( \BB{X}_t \in B,N_{]t_{N(u)+1}, t]} = n \mid \C{F}_{t_{N(u) + 1}}) G_u(dt_{N(u) + 1}, dx_{N(u)+1})\\
&= \int_u^t  \int_{X \setminus \{\BB{X}_u\}} q(\BB{X}_u, s, dy) e^{-\int_u^s q(\BB{X}_u, \theta) d\theta} \bar{P}^{(n)} (s,y;t,B)ds  = \bar{P}^{(n+1)} (u,\BB{X}_u;t,B),
\end{aligned}
\end{multline}
where the first equality holds since $N_{]u,t]} = 1+ N_{]t_{N(u)+1}, t]}$ for $N_{]u,t]} \ge 1$ and
since $\BB{E}(\BB{E}(Z \mid \mathfrak{D})) = \BB{E}(Z)$ for any random variable Z and any $\sigma$-field $\mathfrak{D}$,
the second equality holds since $\sigma(\C{F}_u, t_{N(u) + 1}, x_{N(u)+1}) = \C{F}_{t_{N(u) + 1}}$, the third
equality follows from \eqref{con3a} and \eqref{iii6}, and the last
equality is \eqref{bn}.  Equality \eqref{iii6} is
proved.

Observe that for $u,t \in \BB{R}_+$, $u < t$, $B \in \B(X)$,
\begin{multline}
\label{i5}
\mathbb{P}(\mathbb{X}_t \in B \mid \C{F}_u) = \BB{P}(\BB{X}_t \in B \mid \C{F}_u)I\{u < t_\infty\} + \BB{P}(\BB{X}_t \in B \mid \C{F}_u)I\{u \ge t_\infty\} \\
\begin{aligned}
&= \sum_{n \geq 0} \mathbb{P}(\mathbb{X}_t \in B, N_{]u,t]}= n \mid \C{F}_u)I\{u < t_\infty\} = \sum_{n \geq 0} \bar{P}^{(n)} (u,\BB{X}_u;t,B)I\{u < t_\infty\}\\
&= \bar{P}(u,\BB{X}_u;t,B)I\{u < t_\infty\} = \bar{P}(u,\BB{X}_u;t,B)I\{\BB{X}_u \in X\},
\end{aligned}
\end{multline}
where the first equality holds since $\{\{u < t_\infty\}, \{u \ge t_\infty\}\}$ is a partition of $\Omega$ and $\{u < t_\infty\},\,\{u \ge t_\infty\}  \in \C{F}_u$, the second equality holds since $\BB{X}_t\in X$ implies
$t<t_\infty$, the third equality follows from \eqref{iii6}, the fourth equality follows from \eqref{def}, and
the last one holds since $\{u < t_\infty\} = \{\BB{X}_u \in X\}$. As follows from
\eqref{i5}, the function  $\BB{P}(\mathbb{X}_t \in B \mid
\C{F}_u)$ is $\sigma(\BB{X}_u)$-measurable. Thus,
\begin{equation}
\label{e:MarkProp}
\BB{P}(\mathbb{X}_t \in B \mid \C{F}_u) = \BB{P}(\BB{P}(\mathbb{X}_t \in B \mid \C{F}_u) \mid \BB{X}_u) = \BB{P}(\mathbb{X}_t \in B  \mid \BB{X}_u),
\end{equation}
where the second equality holds because $\sigma(\BB{X}_u)$
$\subseteq$  $\C{F}_u$; see e.g. Br\'emaud~\cite[p. 280]{Bremaud}.
Thus, \eqref{e:MarkProp11} follows from \eqref{i5} and \eqref{e:MarkProp}.
\end{proof}
\end{theorem}

\section{Backward Kolmogorov equation}
\label{S-BKE}
In this section, we show that the transition function
$\bar{P}$ defined in \eqref{def} is the minimal
non-negative solution to the backward Kolmogorov equation.
For a continuous $Q$-function $q$, relevant results were established
by Feller \cite[Theorems 2, 3]{Fel}.

\begin{theorem}
\label{BKE} The function $\bar{P}(u,x;t,B)$ satisfies the following properties:\\
(i) $\bar{P}(u,x;t,B)$ is for fixed $x,t,B$ an absolutely continuous function in $u$ and satisfies uniformly in $B \in \mathfrak{B}(X)$ the
boundary condition
\begin{equation}
\label{CC2} \lim\limits_{u \to t^-}\bar{P}(u,x;t,B) = I \{ x \in B \}.
\end{equation}
(ii) For all $x,t,B$, the function $\bar{P}(u,x;t,B)$
satisfies  for almost every $u<t$ the backward Kolmogorov equation
\begin{equation}
\label{BKDE} \frac{\partial}{\partial u}{P}(u,x;t,B) =
q(x,u){P}(u,x;t,B) - \int_{X \setminus
\{x\}}q(x,u,dy)P(u,y;t,B).
\end{equation}
\begin{proof}
(i) For all $x \in X$, $u,t \in \BB{R}_+$, $u < t$, and $B \in \B(X)$,
\begin{multline}
\label{BKDE-g}
\bar{P}(u,x;t,B) = \sum\limits_{n =0}^{\infty}\bar{P}^{(n)}(u,x;t,B)\\
\begin{aligned}
&= I\{x \in B\}e^{-\int_u^t q(x,s)ds} + \sum_{n = 1}^\infty \int_u^t e^{-\int_u^s q(x,\theta)d\theta}  \int_{X\setminus \{x\}} q(x,s,dy)\bar{P}^{(n-1)}(s, y; t, B)ds \\
&= I\{x \in B\}e^{-\int_u^t q(x,s)ds} +  \int_u^t e^{-\int_u^s q(x,\theta)d\theta}    \int_{X\setminus \{x\}}  q(x,s,dy) \sum_{n = 1}^\infty \bar{P}^{(n-1)}(s, y; t, B)ds \\
&= I\{x \in B\}e^{-\int_u^t q(x,s)ds} +  \int_u^t e^{-\int_u^s
q(x,\theta)d\theta}  \int_{X\setminus \{x\}} q(x,s,dy)\bar{P}(s,
y; t, B)ds,
\end{aligned}
\end{multline}
where the first equality is \eqref{def}, the second equality
follows from \eqref{b0} and \eqref{bn}, the third equality is
obtained by interchanging the integral and sum,  and the last
one follows from \eqref{def}. For fixed $x,t,B$, equation
\eqref{BKDE-g} implies that $\bar{P}(u,x;t,B)$ is the
sum of two absolutely continuous functions in $u$. Thus,
$\bar{P}(u,x;t,B)$ is for fixed $x,t,B$ an absolutely continuous function in $u$.

Observe that $\bar{P}^{(n)}(u,x;t,B) \le \bar{P}(u,x;t,B) \le 1$ for all $n \ge 0$, $x \in X, u,t \in \BB{R}_+$, $u <t$, and $B \in \B(X)$. Then from \eqref{bn},
\begin{equation}
\label{bound}
\bar{P}^{(n)}(u,x;t,B) \leq \int_{u}^{t}e^{-\int_u^s q(x,\sigma)d\sigma}q(x,s)ds, \qquad n \ge 1.
\end{equation}
This inequality and \eqref{b0} imply that, for any stable
$Q$-function $q$,
\begin{equation}
\label{e31}
\lim_{u \to t^-}\bar{P}^{(n)}(u,x;t,B) = 0 \quad
\mbox{ for all } n \ge 1 \quad \mbox{ and } \quad \lim_{u \to
t^-}\bar{P}^{(0)}(u,x;t,B) = I \{x \in B\}
\end{equation}
uniformly with respect to $B$. Thus,  \eqref{def} and \eqref{e31} imply \eqref{CC2}. \\

\noindent
(ii) Since an absolutely continuous real-valued function is
differentiable almost everywhere on its domain, for all $x,t,B$
the function $\bar{P}(u,x;t,B)$ is differentiable in $u$ almost
everywhere on $]0,t[$. By differentiating
\eqref{BKDE-g}, for  almost every $u < t,$
\begin{multline}
\label{BKDE-df}
\begin{aligned}
\frac{\partial}{\partial u }\bar{P}(u,x; t, B) &= I\{x \in B\} e^{-\int_u^t q(x,s)ds} q(x,u) - \int_{X\setminus \{x\}} q(x,u,dy)\bar{P}(u, y; t, B)\\
& \hspace{1in} + \int_u^t \frac{\partial}{\partial u } e^{-\int_u^s q(x,\theta)d\theta}  \int_{X\setminus \{x\}} q(x,s,dy)\bar{P}(s, y; t, B)ds \\
&= I\{x \in B\} e^{-\int_u^t q(x,s)ds} q(x,u) -\int_{X\setminus \{x\}} q(x,u,dy)\bar{P}(u, y; t, B) \\
& \hspace{1in} + \int_u^t  e^{-\int_u^s q(x,\theta)d\theta}  q(x,u) \int_{X\setminus \{x\}} q(x,s,dy)\bar{P}(s, y; t, B)ds.
\end{aligned}
\end{multline}
In view of \eqref{BKDE-g}, the sum of the first and the last terms in
the last expression of \eqref{BKDE-df} is equal to the first term on the right-hand
side of \eqref{BKDE}.
%
\end{proof}
\end{theorem}

As shown in Feller~\cite[Theorem 2]{Fel}, for a stable continuous
$Q$-function $q$, the transition function $\bar{P}$   satisfies
the backward Kolmogorov equation for all $u$, while
Theorem~\ref{BKE}(ii) states that this equation holds for almost
every $u.$ This difference in formulations takes place because the
continuity of the $Q$-function $q$ and the finiteness of each
integrand in the last expression of \eqref{BKDE-g} guarantee the
existence of the derivative $\frac{\partial}{\partial
u}\bar{P}(u,x;t,B)$ for all $u$.

\begin{definition}
A function $P$ with the same domain as $\bar{P}$ is a solution of the backward Kolmogorov equation \eqref{BKDE} if the function $P$ satisfies the properties stated in Theorem \ref{BKE}.
\end{definition}
The next theorem describes the minimal and uniqueness properties
of the solution $\bar{P}$ of the backward Kolmogorov equation
\eqref{BKDE}.
\begin{theorem}
\label{unique-b} The function $\bar{P}$ is the minimal
non-negative solution of the backward Kolmogorov
equation~\eqref{BKDE}. Also, if $\bar{P}$ is a regular transition function (that is, $\bar{P}(u,x;t,X) $ $=1$ for all
$u,x,t$ in the domain of $\bar P$), then $\bar{P}$ is the unique
non-negative solution of the backward Kolmogorov equation
\eqref{BKDE} that is a measure on $(X, \B(X))$ for fixed $u,x,t$
with $u <t$ and takes values in $[0,1]$.
\begin{proof}  The proof of minimality is similar to the proof of Theorem 3 in
Feller~\cite{Fel}.  We provide it here for completeness. Let $P^\ast$ with the same domain as $\bar{P}$ be a
non-negative solution of the
backward Kolmogorov equation \eqref{BKDE}. Integrating
\eqref{BKDE} from $u$ to $t$ and by using the boundary condition \eqref{CC2},
\begin{equation}
\label{BKDE-d}
P^\ast(u,x; t, B) = I\{x \in B\}e^{- \int_{u}^{t}q(x, s) ds} + \int_{u}^{t} \int_{X \setminus \{x\}} e^{ -\int_{u}^{s}
q(x,\theta) d\theta }  q(x,s,dy) P^\ast(s, y; t, B) ds.
\end{equation}
Since the last term of \eqref{BKDE-d} is non-negative,
\begin{equation}
\label{i6}
P^\ast(u, x; t, B) \geq I\{x \in B\}e^{- \int_{u}^{t}q(x, s) ds} = \bar{P}^{(0)}(u,x;t,B),
\end{equation}
where the last equality is \eqref{b0}. For all $u,x,t, B$ with $u < t$,
assume $P^\ast(u, x; t, B)$ $\ge$ $\sum\limits_{m=0}^{n}\bar{P}^{(m)} (u, x; t, B)$ for some $n \geq 0$. Then from \eqref{BKDE-d}
\begin{multline*}
\begin{aligned}
P^\ast(u,x; t, B) &\ge I\{x \in B\}e^{- \int_{u}^{t}q(x, s) ds} + \int_{u}^{t} \int_{X \setminus \{x\}} e^{ -\int_{u}^{s}
q(x,\theta) d\theta }  q(x,s,dy) \sum_{m=0}^{n}\bar{P}^{(m)} (s, y; t, B) ds\\
&= \bar{P}^{(0)}(u,x;t,B) + \sum_{m=0}^{n} \bar{P}^{(m+1)} (u, x; t, B) = \sum_{m=0}^{n+1} \bar{P}^{(m)} (u, x; t, B),
\end{aligned}
\end{multline*}
where the first equality follows from the assumption that $P^\ast(u, x; t, B)$ $\ge$
$\sum\limits_{m=0}^{n}\bar{P}^{(m)} (u, x; t, B)$
for all $u,x,t,B$ with $u < t$, the second equality follows from \eqref{b0} and \eqref{bn}, and the third equality is straightforward. Thus, by induction, $P^\ast(u, x; t, B)$ $\ge$
$\sum\limits_{m=0}^{n}\bar{P}^{(m)} (u, x; t, B)$ for all $n \ge 0$, $ x\in X,$ $u,t \in \BB{R}_+$, $u < t,$ and $B \in \B(X)$, which implies that
$P^\ast(u, x; t, B) \ge \bar{P}(u,x;t,B)$ for all $u,x,t, B$.

To prove the second part of the theorem, let the solution $P^*$ be
a measure on $(X, \B(X))$ for fixed $u,x,t$ and with values in $[0,1]$. Assume that
$P^\ast(u, x; t, B)\neq \bar{P}(u,x;t,B)$ for at least one tuple
$(u,x,t,B)$.  Then,
\[
\begin{split}
P^\ast(u, x; t, X) & = P^\ast(u, x; t, B) + P^\ast(u, x; t, B^c) \\
& > \bar{P}(u, x; t, B) + \bar{P}(u, x; t, B^c) = \bar{P}(u, x; t, X) = 1,
\end{split}
\]
where the  inequality holds because $P^*(u,x,t,\cdot)\ge {\bar
P}(u,x,t,\cdot)$ for all $u,x,t$. Since $P^\ast$ takes values in $[0,1]$,
the assumption that $P^\ast(u, x; t, B)$ $\neq$ $\bar{P}(u, x;
t, B)$ for atleast one tuple $(u,x,t,B)$ leads to a contradiction.
\end{proof}
\end{theorem}

\section{Forward Kolmogorov equation}
\label{S-FKE}

For the forward Kolmogorov equation, this section provides the
results similar to the results on backward Kolmogorov equation in
Section~\ref{S-BKE}.

\begin{theorem}
\label{FKE}
The function $\bar{P}(u,x;t,B)$ satisfies the following properties:\\
(i) $\bar{P}(u,x;t,B)$ is for fixed $u,x,B$ an absolutely continuous function in $t$ and satisfies uniformly in $B \in \mathfrak{B}(X)$ the
boundary condition
\begin{equation}
\label{CC} \lim\limits_{t \to u^+}\bar{P}(u,x;t,B) = I \{ x \in B \}.
\end{equation}
(ii) For all $u,x$, and $q$-bounded sets $B$, the function
$\bar{P}(u,x;t,B)$ satisfies for almost every $t>u$ the forward Kolmogorov
equation
\begin{equation}
\label{FKDE} \frac{\partial}{\partial t}P(u,x; t, B) = - \int_B
q(y,t)P(u,x;t,dy) + \int_{X} q(y,t, B \setminus \{y\})  P(u, x; t,
dy).
\end{equation}
\begin{proof}
(i) For all $x \in X, u,t \in \mathbb{R}_+$, $u < t$, and $B \in
\mathfrak{B}(X)$, equation \eqref{BKDE-g} implies that the
function $\bar{P}(u,x;t,B)$ is absolutely continuous in $t$ for
fixed $u,x,B$. Also, equations \eqref{b0} and
\eqref{bound} imply that, for any stable $Q$-function $q$,
\begin{equation}
\label{e311}
\lim_{t \to u^+}\bar{P}^{(n)}(u,x;t,B) = 0 \quad
\mbox{ for all } n \ge 1 \quad \mbox{ and } \quad \lim_{t \to
u^+}\bar{P}^{(0)}(u,x;t,B) = I \{x \in B\}
\end{equation}
uniformly with respect to $B$. Thus,  \eqref{def} and \eqref{e311} imply \eqref{CC}. \\

\noindent
(ii) Consider the following non-negative function defined on the domain of $\bar{P}$ by
\begin{equation}
\label{pi}
\Pi(u,x;t,B) = \int_{B \setminus \{x\}} q(x,u,dy) e^{- \int_u^t q(y,s)ds}.
\end{equation}
According to  Feller \cite[Theorem 4]{Fel},  the function
$\bar{P}^{(n)}(u,x; t, B),$  $n\ge 1$, satisfies the recursion
\begin{equation}
\label{en} {\bar P}^{(n)}(u,x;t,B) = \int_{u}^{t}  \int_X \Pi(s,y;t,B){\bar
P}^{(n-1)}(u, x; s, dy) ds.
\end{equation}
Though the function $\bar{P}^{(n)}(u,x; t, B)$ is defined for
continuous $Q$-functions in Feller \cite{Fel}, the proof given
there is correct for Borel $Q$-functions. From \eqref{b0},
\eqref{def}, and \eqref{en},
\begin{equation}
\label{FKDE-g2}
\begin{aligned}
\bar{P}(u,x;t,B)&= \sum_{n = 0}^{\infty}\bar{P}^{(n)}(u,x; t, B) \\
&=  I\{x \in B\}e^{- \int_{u}^{t}q(x, s) ds} + \sum_{n=1}^\infty \int_{u}^{t}  \int_X \Pi (s,y; t, B) \bar{P}^{(n-1)}(u,x; s, dy) ds \\
&=  I\{x \in B\}e^{- \int_{u}^{t}q(x, s) ds} + \int_{u}^{t}  \int_X \Pi (s,y; t, B) \sum_{n = 1}^{\infty}\bar{P}^{(n-1)}(u,x; s, dy) ds \\
&= I\{x \in B\}e^{- \int_{u}^{t}q(x, s) ds} + \int_{u}^{t} \int_X
\Pi (s,y; t, B) \bar{P}(u,x; s, dy) ds.
\end{aligned}
\end{equation}
 Since
$\bar{P}(u,x;t,B)$ is an absolutely continuous function in $t$
for fixed $u,x,B$, the derivative  $\frac{\partial}{\partial t}\bar{P}(u,x; t, B)$
exists for almost every $t \in \, ]u,\infty[$.
By differentiating \eqref{FKDE-g2},
for almost every $t>u$,
\begin{multline}
\label{df2}
\frac{\partial}{\partial t}\bar{P}(u,x; t, B) = - I\{x
\in B\}  e^{- \int_{u}^{t}q(x, s) ds} q(x,t) \\
+  \int_X \Pi(t,y; t, B) \bar{P}(u,x;t,dy) + \int_{u}^{t}\frac{\partial}{\partial t} \int_X \Pi(s,y; t, B) \bar{P}(u,x;s,dy) ds.
\end{multline}
By differentiating \eqref{pi} with respect to $t$, for all
$q$-bounded sets $B \in \mathcal{B}(X)$,
\begin{equation}
\label{pi-df}
\frac{\partial}{\partial t}\Pi(u,x;t,B)= \int_{B \setminus \{x\}} q(x,u,dy)  \frac{\partial}{\partial t}
e^{- \int_u^t q(y,s)ds} = - \int_{B} q(y,t)
\Pi(u,x;t,dy).
\end{equation}
Combining \eqref{df2} and \eqref{pi-df} and observing that $\Pi(t,y; t, B) = q(y,t, B\setminus \{y\})$,  for all $q$-bounded sets $B$,
\begin{multline}
\label{df3}
\frac{\partial}{\partial t}\bar{P}(u,x; t, B) = - I\{x \in B\} e^{- \int_{u}^{t}q(x, s) ds } q(x,t)\\
+ \int_X q(y,t, B\setminus \{y\}) \bar{P}(u,x;t,dy) -
\int_{u}^{t}  \int_X \int_{B} q(z,t) \Pi(s,y;t,dz)
\bar{P}(u,x;s,dy) ds,
\end{multline}
for almost every $t > u$. By substituting $\bar{P}(u,x;t,dz)$ in
the left-hand side of  the following equality with the final expression in \eqref{FKDE-g2},
\begin{equation}
\begin{aligned}
\label{df4}
\int_B q(z,t)\bar{P}(u,x;t,dz) &= I\{x \in B\}e^{- \int_{u}^{t}q(x, s) ds } q(x,t) \\
&+ \int_u^t \int_X \int_B q(z,t) \Pi(s,y;t,dz) \bar{P}(u,x;s,dy) ds.
\end{aligned}
\end{equation}
Formulae \eqref{df3} and \eqref{df4} imply statement (ii) of the theorem.
%
\end{proof}
\end{theorem}

\begin{corollary}[Feller~{\cite[Equation (37)]{Fel}}]
\label{Cor}
For all $x \in X, u, t \in \mathbb{R}_+$, $u < t$, and $q$-bounded sets $B \in \mathfrak{B}(X)$, the function $\bar{P}(u,x;t,B)$ defined in \eqref{def} satisfies
\begin{multline}
\label{FKDE-g}
P(u,x;t,B) = I\{x \in B\} \\
+ \int_{u}^{t} ds \int_X q(y,s,B \setminus \{y\}) P(u,x; s, dy) -\int_{u}^{t} ds \int_{B} q(y,s) P(u,x;s,dy).
\end{multline}
\begin{proof}
Integrating \eqref{FKDE} from $u$ to $t$ and by using the boundary condition \eqref{CC}, we get \eqref{FKDE-g}. Thus, by Theorem \ref{FKE}, for all $x \in X, u, t \in \mathbb{R}_+$, $u < t$, and $q$-bounded sets $B \in \mathfrak{B}(X)$, the function $\bar{P}(u,x;t,B)$ satisfies \eqref{FKDE-g}.
\end{proof}
\end{corollary}

\begin{definition}
A function $P$ with the same domain as $\bar{P}$ is a solution of the forward Kolmogorov equation \eqref{FKDE} if the function $P$ satisfies the properties stated in Theorem \ref{FKE}.
\end{definition}
Following the proof of Theorem \ref{unique-b}, we establish the minimal and uniqueness properties of the solution $\bar{P}$ of the forward Kolmogorov equation \eqref{FKDE} in theorem \ref{unique-f}. We remark that the function $\bar{P}$ is the minimal non-negative solution of \eqref{FKDE} on a restricted domain $(u,t \in \BB{R}_+, u <t, x \in X,$ and $q$-bounded sets $B \in \B(X))$.

\begin{theorem}
\label{unique-f}
The function $\bar{P}(u,x;t,B)$, being restricted to $q$-bounded sets $B$,
is the minimal non-negative solution of
the forward Kolmogorov equation \eqref{FKDE}. Also, if $\bar{P}$ is a
regular transition function (that is, $\bar{P}(u,x;t,X)=1$ for
all $u,x,t$ in the domain of $\bar P$), then $\bar{P}$ is the
unique non-negative solution of the forward Kolmogorov equation \eqref{FKDE}
that is a measure on $(X, \B(X))$ for fixed $u,x,t$ with $u < t$ and takes values in $[0,1]$.
\begin{proof}
Let $P^\ast$ defined on the same domain as $\bar{P}$ be a non-negative solution of the forward Kolmogorov equation
\eqref{FKDE}. Integrating \eqref{FKDE} from $u$ to $t$ and
by using the boundary condition \eqref{CC}, for all $x \in X$, $u, t\in \BB{R}_+$ with  $u < t$, and $q$-bounded sets $B \in \B(X)$,
\begin{equation}
\label{FKDE-d}
P^\ast(u,x; t, B) = I\{x \in B\}e^{- \int_{u}^{t}q(x, s) ds} + \int_{u}^{t} ds \int_X \Pi (s,y; t, B) P^\ast(u,x; s, dy),
\end{equation}
for all $q$-bounded sets $B$. Since the last term of \eqref{FKDE-d} is non-negative,
\begin{equation}
\label{i6}
P^\ast(u, x; t, B) \geq I\{x \in B\}e^{- \int_{u}^{t}q(x, s) ds} = \bar{P}^{(0)}(u,x;t,B),
\end{equation}
where the last equality is \eqref{b0}. For all $x,u,t$ with $u < t$ and $q$-bounded sets $B$,
assume $P^\ast(u, x; t, B)$ $\ge \sum\limits_{m=0}^{n}\bar{P}^{(m)}(u,x;t,B)$ for some $n \ge 0$. Then from \eqref{FKDE-d}
\begin{multline*}
\begin{aligned}
P^\ast(u,x; t, B) &\ge I\{x \in B\}e^{- \int_{u}^{t}q(x, s) ds} + \int_{u}^{t} ds \int_X \Pi (s,y; t, B) \sum_{m=0}^{n}\bar{P}^{(m)}(u,x;s,dy) \\
&= \bar{P}^{(0)}(u,x;t,B) + \sum_{m=0}^{n}\bar{P}^{(m+1)}(u,x;t,B) = \sum_{m=0}^{n+1}\bar{P}^{(m)}(u,x;t,B).
\end{aligned}
\end{multline*}
Thus, by induction, $P^\ast(u,
x; t, B) \ge \sum_{m =0}^{n}\bar{P}^{(m)}(u,x;t,B)$ for all $n \ge 0$, $x \in X$, $u, t\in \BB{R}_+$ with $u < t$,
and $q$-bounded sets $B \in \B(X)$, which
implies that $P^\ast(u, x; t, B) \ge \bar{P}(u,x;t,B)$ for all
$u,x,t$ with $u <t$, and for all $q$-bounded sets $B$.

To prove the uniqueness property of $\bar{P}$, let the solution $P^\ast$ be a measure on $(X, \B(X))$
for fixed $u,x,t$ with $u<t$ and with values in $[0,1]$. It follows from statement (i)
of the theorem that for all $B \in
\mathfrak{B}(X)$
\begin{equation}
\label{m-B}
P^\ast(u, x; t, B) = \lim\limits_{n \to \infty}P^\ast(u, x; t, B
\cap B_n) \geq \lim\limits_{n \to \infty}\bar{P}(u, x; t, B \cap
B_n) = \bar{P}(u, x; t, B),
\end{equation}
where $\{B_n\}_{n \geq 0}$ is  an increasing sequence of
$q$-bounded sets  such that $B_n \uparrow X$ as $n \uparrow
\infty$, whose existence is guaranteed by Assumption \ref{A}.
If $\bar{P}(u,x;t,X) = 1$ for all $u,x,t$, then the uniqueness of $\bar{P}$ within the set of
solutions to the forward Kolmogorov equation that take values in $[0,1]$ and
that are measures on $(X, \B(X))$ for fixed $u,x,t$ with $u <t$ follows from
the minimality of $\bar{P}$ \eqref{m-B} and from the same arguments as in the
proof of uniqueness in Theorem \ref{unique-b}.
\end{proof}
\end{theorem}

\begin{remark}\label{rm4.1} {\rm The results of this paper can be extended to
non-conservative $Q$-functions.  As mentioned in
section~\ref{sec1}, any non-conservative $Q$-function $q$ can be
transformed into a conservative $Q$-function by adding a state
$\tilde{x}$ to $X$ with $q(x,t,\{\tilde{x}\}):= -q(x,t,X)$,
$q(\tilde{x}, t, X):= 0$, and $q(\tilde{x}, t, \{\tilde{x}\}) :=
0$, where $x \in X$ and $t \in \mathbb{R}_+$.  According to
Theorem~\ref{JMP}, there is a transition function $\bar P$ of a
jump Markov process with the state space ${\tilde
X}=X\cup\{\tilde{x}\}$, and this process is determined by the
initial state distribution and by the compensator defined by the
modified $Q$-function.
The proofs of the results of sections~\ref{S-BKE} and \ref{S-FKE}
do not use the assumption that the $Q$-function $q$ is
conservative.  Therefore, these results remain valid for
non-conservative $Q$-functions.  However, the validity of the
condition ${\bar P}(u,x;t,X)=1$ for all $x,u,t$ with $u <t$
in Theorems~\ref{unique-b} and  \ref{unique-f} is possible only if
$q(x,t,X)=0$ almost everywhere in $t$ for each $x\in X$.  Thus, in
fact, $q$ is conservative, if $\bar{P}(u,x;t,X) = 1$ for all $x,u,t$ with $u<t$.
It is also easy to see that the minimal non-negative
solutions of both the backward and forward Kolmogorov equations
are equal to ${\bar P}(u,x;t,B),$ when $x\in X$ and $B\in \B(X)$,
where the transition function $\bar P$ is described in the
previous paragraph for a broader domain.} \end{remark}

\begin{remark}\label{rm4.2} {\rm  In this paper, transition
functions $P(u,x;t,B)$ and $\bar P(u,x;t,B)$ are defined for
$u>0.$ All the results of Sections 3 and 4 hold for $u\ge 0$ with
the same proofs. When $x$ is the initial state of the process
$\BB{X}$, the results of Section 2 also hold for $u\ge 0$.}
\end{remark}

\bibliographystyle{elsarticle-num}

\end{document}